\documentclass{amsart}
\usepackage[latin1]{inputenc}

\usepackage{amssymb,amsmath}
\usepackage{amsthm}
\usepackage{paralist}
\usepackage{url}

\newtheorem{theorem}{Theorem}

\newtheorem{lemma}[theorem]{Lemma}

\newcommand{\minor}[3]{\ensuremath{{#1}_{{#2} \gets {#3}}}} 
\DeclareMathOperator{\ess}{ess} 
\DeclareMathOperator{\essl}{ess^{<}} 
\DeclareMathOperator{\gap}{gap} 

\begin{document}
\title[The effect of variable identification on essential arity]{On the effect of variable identification on the essential arity of functions}
\author{Miguel Couceiro}
\address{Department of Mathematics, Statistics and Philosophy, University of Tampere, FI-33014 Tampereen yliopisto, Finland}
\email{miguel.couceiro@uta.fi}
\author{Erkko Lehtonen}
\address{Institute of Mathematics, Tampere University of Technology, P.O. Box 553, FI-33101 Tampere, Finland}
\email{erkko.lehtonen@tut.fi}
\date{\today}
\begin{abstract}
We show that every function of several variables on a finite set of $k$ elements with $n > k$ essential variables has a variable identification minor with at least $n - k$ essential variables. This is a generalization of a theorem of Salomaa on the essential variables of Boolean functions. We also strengthen Salomaa's theorem by characterizing all the Boolean functions $f$ having a variable identification minor that has just one essential variable less than $f$.
\end{abstract}
\maketitle

\section{Introduction}

Theory of essential variables of functions has been developed by several authors \cite{Chimev,EKMM,Salomaa,Yablonski}. In this paper, we discuss the problem how the number of essential variables is affected by identification of variables (diagonalization). Salomaa \cite{Salomaa} proved the following two theorems: one deals with operations on arbitrary finite sets, while the other deals specifically with Boolean functions. We denote the number of essential variables of $f$ by $\ess f$.

\begin{theorem}
Let $A$ be a finite set with $k$ elements. For every $n \leq k$, there exists an $n$-ary operation $f$ on $A$ such that $\ess f = n$ and every identification of variables produces a constant function.
\end{theorem}

Thus, in general, essential variables can be preserved when variables are identified only in the case that $n > k$.

\begin{theorem}
\label{thm:SalomaaMain}
For every Boolean function $f$ with $\ess f \geq 2$, there is a function $g$ obtained from $f$ by identification of variables such that $\ess g \geq \ess f - 2$.
\end{theorem}

Identification of variables together with permutation of variables and cylindrification induces a quasi-order on operations whose relevance has been made apparent by several authors \cite{Couceiro,EFHH,FH,Lehtonen,Pippenger,Wang,Zverovich}. In the case of Boolean functions, this quasi-order was studied in \cite{CP} where Theorem \ref{thm:SalomaaMain} was fundamental in deriving certain bounds on the essential arity of functions.

In this paper, we will generalize Theorem \ref{thm:SalomaaMain} to operations on arbitrary finite sets in Theorem \ref{thm:gen}. We will also strengthen Theorem \ref{thm:SalomaaMain} on Boolean functions in Theorem \ref{thm:str} by determining the Boolean functions $f$ for which there exists a function $g$ obtained from $f$ by identification of variables such that $\ess g = \ess f - 1$.

\section{Variable identification minors}

Let $A$ and $B$ be arbitrary nonempty sets. A \emph{$B$-valued function of several variables on $A$} is a mapping $f : A^n \to B$ for some positive integer $n$, called the \emph{arity} of $f$. $A$-valued functions on $A$ are called \emph{operations on $A$.} Operations on $\{0,1\}$ are called \emph{Boolean functions.}

We say that the $i$-th variable is \emph{essential} in $f$, or $f$ depends on $x_i$, if there are elements $a_1, \ldots, a_n, b \in A$ such that
\[
f(a_1, \ldots, a_i, \ldots, a_n) \neq f(a_1, \ldots, a_{i-1}, b, a_{i+1}, \ldots, a_n).
\]
The number of essential variables in $f$ is called the \emph{essential arity} of $f$, and it is denoted by $\ess f$. Thus the only functions with essential arity zero are the constant functions.

For an $n$-ary function $f$, we say that an $m$-ary function $g$ is obtained from $f$ by \emph{simple variable substitution} if there is a mapping $\sigma : \{1, \ldots, n\} \to \{1, \ldots, m\}$ such that
\[
g(x_1, \ldots, x_m) = f(x_{\sigma(1)}, \ldots, x_{\sigma(n)}).
\]
In the particular case that $n = m$ and $\sigma$ is a permutation of $\{1, \ldots, n\}$, we say that $g$ is obtained from $f$ by \emph{permutation of variables.} For indices $i, j \in \{1, \ldots, n\}$, $i \neq j$, if $x_i$ and $x_j$ are essential in $f$, then the function $\minor{f}{i}{j}$ obtained from $f$ by the simple variable substitution
\[
\minor{f}{i}{j}(x_1, \ldots, x_n) = f(x_1, \ldots, x_{i-1}, x_j, x_{i+1}, \ldots, x_n)
\]
is called a \emph{variable identification minor} of $f$, obtained by identifying $x_i$ with $x_j$. Note that $\ess \minor{f}{i}{j} < \ess f$, because $x_i$ is not essential in $\minor{f}{i}{j}$ even though it is essential in $f$.

We define a quasiorder on the set of all $B$-valued functions of several variables on $A$ as follows: $f \leq g$ if and only if $f$ is obtained from $g$ by simple variable substitution. If $f \leq g$ and $g \leq f$, we denote $f \equiv g$. If $f \leq g$ but $g \not\leq f$, we denote $f < g$. It can be easily observed that if $f \leq g$ then $\ess f \leq \ess g$, with equality if and only if $f \equiv g$.

For a $B$-valued function $f$ of several variables on $A$, we denote the maximum essential arity of a variable identification minor of $f$ by
\[
\essl f = \max_{g < f} \ess g,
\]
and we define the \emph{arity gap} of $f$ by $\gap f = \ess f - \essl f$.

\section{Generalization of Theorem \ref{thm:SalomaaMain}}

\begin{theorem}
\label{thm:gen}
Let $A$ be a finite set of $k \geq 2$ elements, and let $B$ be a set with at least two elements. Every $B$-valued function of several variables on $A$ with $n > k$ essential variables has a variable identification minor with at least $n - k$ essential variables.
\end{theorem}

In the proof of Theorem \ref{thm:gen}, we will make use of the following theorem due to Salomaa \cite[Theorem 1]{Salomaa}, which is a strengthening of Yablonski's \cite{Yablonski} ``fundamental lemma''.
\begin{theorem}
\label{thm:SalomaaAux}
Let the function $f : M_1 \times \cdots \times M_n \to N$ depend essentially on all of its $n$ variables, $n \geq 2$. Then there is an index $j$ and an element $c \in M_j$ such that the function
\[
f(x_1, \ldots, x_{j-1}, c, x_{j+1}, \ldots, x_n)
\]
depends essentially on all of its $n - 1$ variables.
\end{theorem}

We also need the following auxiliary lemma.
\begin{lemma}
\label{lem}
Let $f$ be an $n$-ary function with $\ess f = n > k$. Then there are indices $1 \leq i < j \leq k+1$ such that at least one of the variables $x_1$, \ldots, $x_{k+1}$ is essential in $\minor{f}{i}{j}$.
\end{lemma}
\begin{proof}
Since $x_1$ is essential in $f$, there are elements $a_1, \ldots, a_n, b \in A$ such that
\[
f(a_1, a_2, \ldots, a_n) \neq f(b, a_2, \ldots, a_n).
\]
Thus there are indices $1 \leq i < j \leq k+1$ such that $a_i = a_j$. If $i \neq 1$, then it is clear that $x_1$ is essential in $\minor{f}{i}{j}$. If there are no such $i$ and $j$ with $i \neq 1$, then $i = 1 < j$ and we have that $b = a_l$ for some $1 < l \leq k+1$, $l \neq j$. For $m = 1, \ldots, n$, let $c_m = a_m$ if $m \notin \{1, j, l\}$ and let $c_m = a_1$ if $m \in \{1, j, l\}$. Then $f(c_1, c_2, \ldots, c_n)$ is distinct from at least one of $f(a_1, a_2, \ldots, a_n)$ and $f(b, a_2, \ldots, a_n)$. If $f(c_1, c_2, \ldots, c_n) \neq f(a_1, a_2, \ldots, a_n)$, then $x_l$ is essential in $\minor{f}{1}{j}$. If $f(c_1, c_2, \ldots, c_n) \neq f(b, a_2, \ldots, a_n)$, then $x_l$ is essential in $\minor{f}{1}{l}$.
\end{proof}

\begin{proof}[Proof of Theorem \ref{thm:gen}]
By Theorem \ref{thm:SalomaaAux}, there exist $k + 1$ constants $c_1, \ldots, c_{k+1} \in A$ such that, after a suitable permutation of variables, the function
\[
f(c_1, \ldots, c_{k+1}, x_{k+2}, \ldots, x_n)
\]
depends on all of its $n - k - 1$ variables. There are indices $1 \leq i < j \leq k+1$ such that $c_i = c_j$, and by Lemma \ref{lem} there are indices $1 \leq l < m \leq k+1$ such that at least one of the variables $x_1$, \ldots, $x_{k+1}$ is essential in $\minor{f}{l}{m}$. With a suitable permutation of variables, we may assume that $i = 1$, $j = 2$, $1 \leq l \leq 3$, $m = l+1$.

If one of the variables $x_1, \ldots, x_{k+1}$ is essential in $\minor{f}{1}{2}$, then we are done. Otherwise we have that for all $a_{k+2}, \ldots, a_n \in A$,
\[
f(c_1, c_1, c_3, c_4, \ldots, c_{k+1}, a_{k+2}, \ldots, a_n) = f(c_3, c_3, c_3, c_4, \ldots, c_{k+1}, a_{k+2}, \ldots, a_n).
\]
Thus the variables $x_{k+2}, \ldots, x_n$ are essential in $\minor{f}{2}{3}$. If one of the variables $x_1, \ldots, x_{k+1}$ is essential in $\minor{f}{2}{3}$, then we are done. Otherwise we have that for all $a_{k+2}, \ldots, a_n \in A$,
\[
f(c_3, c_3, c_3, c_4, \ldots, c_{k+1}, a_{k+2}, \ldots, a_n) = f(c_3, c_4, c_4, c_4, \ldots, c_{k+1}, a_{k+2}, \ldots, a_n),
\]
and so the variables $x_{k+2}, \ldots, x_n$ are essential in $\minor{f}{3}{4}$ and also at least one of $x_1, \ldots, x_{k+1}$ is essential in $\minor{f}{3}{4}$.
\end{proof}

We would like to remark that our proof is considerably simpler than Salomaa's original proof of Theorem \ref{thm:SalomaaMain}.

\section{Strengthening of Theorem \ref{thm:SalomaaMain}}

It is well-known that every Boolean function is represented by a unique multilinear polynomial over the two-element field. Such a representation is called the \emph{Zhegalkin polynomial} of $f$. It is clear that a variable is essential in $f$ if and only if it occurs in the Zhegalkin polynomial of $f$. We denote by $\deg \mathfrak{p}$ the degree of polynomial $\mathfrak{p}$. If $\mathfrak{p}$ is the Zhegalkin polynomial of $f$, then we denote the Zhegalkin polynomial of $\minor{f}{i}{j}$ by $\minor{\mathfrak{p}}{i}{j}$. Note that the only polynomials of degree $0$ are the constant polynomials.

\begin{theorem}
\label{thm:str}
Let $f$ be a Boolean function with at least $2$ essential variables. Then the arity gap of $f$ is $2$ if and only if the Zhegalkin polynomial of $f$ is of one of the following special forms:
\begin{itemize}
\item $x_{i_1} + x_{i_2} + \dots + x_{i_n} + c$,
\item $x_i x_j + x_i + c$,
\item $x_i x_j + x_i x_k + x_j x_k + c$,
\item $x_i x_j + x_i x_k + x_j x_k + x_i + x_j + c$,
\end{itemize}
where $c \in \{0,1\}$. Otherwise the arity gap of $f$ is $1$.
\end{theorem}

We prove first an auxiliary lemma that takes care of the functions of essential arity at least $4$ whose Zhegalkin polynomial has degree $2$.

\begin{lemma}
\label{lemma:deg2}
If $f$ is a Boolean function with at least four essential variables and the Zhegalkin polynomial of $f$ has degree two, then the arity gap of $f$ is one.
\end{lemma}
\begin{proof}
Denote the Zhegalkin polynomial of $f$ by $\mathfrak{p}$. We need to consider several cases and subcases.

\emph{Case 1.} Assume first that $\mathfrak{p}$ is of the form
\[
\mathfrak{p} = x_i x_j + x_i x_k + x_j x_k + x_i \mathfrak{a}_i + x_j \mathfrak{a}_j + x_k \mathfrak{a}_k + \mathfrak{a},
\]
where $\mathfrak{a}_i$, $\mathfrak{a}_j$, $\mathfrak{a}_k$ are polynomials of degree at most $1$ and $\mathfrak{a}$ is a polynomial of degree at most $2$ such that there are no occurrences of variables $x_i$, $x_j$, $x_k$ in $\mathfrak{a}_i$, $\mathfrak{a}_j$, $\mathfrak{a}_k$, $\mathfrak{a}$.

\emph{Subcase 1.1.} Assume that $\deg \mathfrak{a}_i = \deg \mathfrak{a}_j = \deg \mathfrak{a}_k = 0$. Then $\mathfrak{a}$ contains a variable $x_l$ distinct from $x_i$, $x_j$, $x_k$, and we can write $\mathfrak{a} = x_l \mathfrak{a}' + \mathfrak{a}''$, where $\mathfrak{a}'$ and $\mathfrak{a}''$ do not contain $x_l$. Then $\minor{f}{l}{i}$ is represented by the polynomial
\[
\minor{\mathfrak{p}}{l}{i} = x_i x_j + x_i x_k + x_j x_k + x_i \mathfrak{a}' + \mathfrak{a}'',
\]
where all essential variables of $f$ except for $x_l$ occur, because no terms cancel, and hence $\gap f = 1$.

\emph{Subcase 1.2.} Assume that at least one of $\mathfrak{a}_i$, $\mathfrak{a}_j$, $\mathfrak{a}_k$ has degree $1$, say $\deg \mathfrak{a}_i = 1$. Then $\mathfrak{a}_i$ contains a variable $x_l$ distinct from $x_i$, $x_j$, $x_k$, and so $\mathfrak{a}_i = x_l + \mathfrak{a}'_i$, where $\mathfrak{a}'_i$ has degree at most $1$ and does not contain $x_l$. Consider
\[
\minor{\mathfrak{p}}{j}{k} = x_k (1 + \mathfrak{a}_j + \mathfrak{a}_k) + x_i \mathfrak{a}_i + \mathfrak{a}.
\]
If all essential variables of $f$ except for $x_j$ occur in $\minor{\mathfrak{p}}{j}{k}$, then $\gap f = 1$ and we are done. Otherwise we need to analyze three different subcases.

\emph{Subcase 1.2.1.} Assume that variable $x_k$ occurs in $\minor{\mathfrak{p}}{j}{k}$ but there is a variable $x_l$ that occurs in $\mathfrak{a_j}$ and $\mathfrak{a}_k$ but not in $\mathfrak{a}_i$ nor in $\mathfrak{a}$ such that $x_l$ does not occur in $\minor{\mathfrak{p}}{j}{k}$ (due to some cancelling terms in $\mathfrak{a}_j$ and $\mathfrak{a}_k$). Write $\mathfrak{a}_j = x_l + \mathfrak{a}'_j$, $\mathfrak{a}_k = x_l + \mathfrak{a}'_k$, and consider
\[
\begin{split}
\minor{\mathfrak{p}}{j}{l} &= x_i x_l + x_i x_k + x_l x_k + x_i \mathfrak{a}_i + x_l + x_l \mathfrak{a}'_j + x_k x_l + x_k \mathfrak{a}'_k + \mathfrak{a} \\
&= x_i x_l + x_i x_k + x_i \mathfrak{a}_i + x_l + x_l \mathfrak{a}'_j + x_k \mathfrak{a}'_k + \mathfrak{a}.
\end{split}
\]
Every essential variable of $f$ except for $x_j$ occurs in $\minor{\mathfrak{p}}{j}{l}$, and hence $\gap f = 1$.

\emph{Subcase 1.2.2.} Assume that $x_k$ does not occur in $\minor{\mathfrak{p}}{j}{k}$. In this case $\mathfrak{a}_j = \mathfrak{a}_k + 1$. Consider
\[
\minor{\mathfrak{p}}{j}{i} = x_i (1 + \mathfrak{a}_i + \mathfrak{a}_j) + x_k \mathfrak{a}_k + \mathfrak{a}.
\]
If any term of $\mathfrak{a}_j$ is cancelled by a term of $\mathfrak{a}_i$, it still remains as a term of $\mathfrak{a}_k$, and hence all variables occurring in $\mathfrak{a}_i$, $\mathfrak{a}_j$, $\mathfrak{a}_k$ occur in $\minor{\mathfrak{p}}{j}{i}$. If both $x_i$ and $x_k$ also occur in $\minor{\mathfrak{p}}{j}{i}$, then all essential variables of $f$ except for $x_j$ occur in $\minor{\mathfrak{p}}{j}{i}$, and so $\gap f = 1$.

If $x_k$ does not occur in $\minor{\mathfrak{p}}{j}{i}$, then $\mathfrak{a}_k = 0$ and so $\mathfrak{a}_j = 1$. Then
\[
\minor{\mathfrak{p}}{l}{i} = x_i x_j + x_i x_k + x_j x_k + x_i + x_i \mathfrak{a}'_i + x_j + \mathfrak{a},
\]
and every essential variable of $f$ except for $x_l$ occurs in $\minor{\mathfrak{p}}{l}{i}$. Thus $\gap f = 1$.

If $x_i$ does not occur in $\minor{\mathfrak{p}}{j}{i}$, then $\mathfrak{a}_j = \mathfrak{a}_i + 1$, and hence $\mathfrak{a}_i = \mathfrak{a}_k$. Consider then
\[
\minor{\mathfrak{p}}{i}{k} = x_k (1 + \mathfrak{a}_i + \mathfrak{a}_k) + x_j \mathfrak{a}_j + \mathfrak{a}
= x_k + x_j \mathfrak{a}_j + \mathfrak{a}.
\]
Again all essential variables of $f$ except for $x_i$ occur in $\minor{\mathfrak{p}}{i}{k}$, and so $\gap f = 1$.

\emph{Subcase 1.2.3.} Assume that both $x_i$ and $x_k$ occur in $\minor{\mathfrak{p}}{j}{k}$ but there is a variable $x_l$ occurring in $\mathfrak{a}_i$ and in $\mathfrak{a}_j$ but not in $\mathfrak{a}_k$ nor in $\mathfrak{a}$ such that $x_l$ does not occur in $\minor{\mathfrak{p}}{j}{k}$ (due to some cancelling terms in $\mathfrak{a}_i$ and $\mathfrak{a}_j$). Write $\mathfrak{a}_i = x_l + \mathfrak{a}'_i$, $\mathfrak{a}_j = x_l + \mathfrak{a}'_j$, and consider
\[
\begin{split}
\minor{\mathfrak{p}}{j}{l} &= x_i x_l + x_i x_k + x_l x_k + x_i x_l + x_i \mathfrak{a}'_i + x_l + x_l \mathfrak{a}'_j + x_k \mathfrak{a}_k + \mathfrak{a} \\
&= x_i x_k + x_l x_k + x_i \mathfrak{a}'_i + x_l + x_l \mathfrak{a}'_j + x_k \mathfrak{a}_k + \mathfrak{a}.
\end{split}
\]
Every essential variable of $f$ except for $x_j$ occurs in $\minor{\mathfrak{p}}{j}{l}$, and so $\gap f = 1$.

\emph{Case 2.} Assume then that $\mathfrak{p}$ is of the form
\[
\mathfrak{p} = x_i x_j + x_i x_k \mathfrak{a}_{ik} + x_i \mathfrak{a}_i + x_j \mathfrak{a}_j + x_k \mathfrak{a}_k + \mathfrak{a}, 
\]
where $\mathfrak{a}_{ik}$ is a polynomial of degree $0$; $\mathfrak{a}_i$, $\mathfrak{a}_j$, $\mathfrak{a}_k$ are polynomials of degree at most $1$; and $\mathfrak{a}$ is a polynomial of degree at most $2$ such that variables $x_i$, $x_j$, $x_k$ do not occur in $\mathfrak{a}_{ik}$, $\mathfrak{a}_i$, $\mathfrak{a}_j$, $\mathfrak{a}_k$, $\mathfrak{a}$.
Note that $\mathfrak{a}_{ik}$ and $\mathfrak{a}_k$ cannot both be $0$, for otherwise $x_k$ would not occur in $\mathfrak{p}$. Consider
\[
\minor{\mathfrak{p}}{j}{i} = x_i (1 + \mathfrak{a}_i + \mathfrak{a}_j) + x_i x_k \mathfrak{a}_{ik} + x_k \mathfrak{a}_k + \mathfrak{a}.
\]
By the above observation that $\mathfrak{a}_{ik}$ and $\mathfrak{a}_k$ are not both $0$, $x_k$ occurs in $\minor{\mathfrak{p}}{j}{i}$. If all essential variables of $f$ except for $x_j$ occur in $\minor{\mathfrak{p}}{j}{i}$, then $\gap f = 1$ and we are done. Otherwise we distinguish between two cases.

\emph{Subcase 2.1.} Assume that $x_i$ does not occur in $\minor{\mathfrak{p}}{j}{i}$. In this case $\mathfrak{a}_j = \mathfrak{a}_i + 1$, $\mathfrak{a}_{ik} = 0$, and $\mathfrak{a}_k \neq 0$. Consider
\[
\begin{split}
\minor{\mathfrak{p}}{i}{k} &= x_j x_k + x_k \mathfrak{a}_{ik} + x_k \mathfrak{a}_i + x_j \mathfrak{a}_j + x_k \mathfrak{a}_k + \mathfrak{a} \\
&= x_j x_k + x_k (\mathfrak{a}_i + \mathfrak{a}_k) + x_j + x_j \mathfrak{a}_i + \mathfrak{a}.
\end{split}
\]
Both $x_j$ and $x_k$ occur in $\minor{\mathfrak{p}}{i}{k}$, because the term $x_j x_k$ cannot be cancelled. If any term of $\mathfrak{a}_i$ is cancelled by a term of $\mathfrak{a}_k$, it still remains in $x_j \mathfrak{a}_i$. Thus, all essential variables of $f$ except for $x_i$ occur in $\minor{\mathfrak{p}}{i}{k}$, and hence $\gap f = 1$.

\emph{Subcase 2.2.} Assume that $x_i$ occurs in $\minor{\mathfrak{p}}{j}{i}$ but there is a variable $x_l$ occurring in $\mathfrak{a}_i$ and $\mathfrak{a}_j$ but not in $\mathfrak{a}_{ik}$, $\mathfrak{a}_k$, nor in $\mathfrak{a}$ such that $x_l$ does not occur in $\minor{\mathfrak{p}}{j}{i}$ (due to some cancelling terms in $\mathfrak{a}_i$ and $\mathfrak{a}_j$). Consider
\[
\minor{\mathfrak{p}}{k}{l} = x_i x_j + x_i x_l \mathfrak{a}_{ik} + x_i \mathfrak{a}_i + x_j \mathfrak{a}_j + x_l \mathfrak{a}_k + \mathfrak{a}.
\]
If $\mathfrak{a}_{ik} = 1$, then the terms $x_i x_l$ in $x_i \mathfrak{a}_i$ and in $x_i x_l \mathfrak{a}_{ik}$ cancel each other. These are the only terms that may be cancelled out. Nevertheless, $x_l$ occurs also in $\mathfrak{a}_j$, and so all essential variables of $f$ except for $x_k$ occur in $\minor{\mathfrak{p}}{k}{l}$. Therefore $\gap f = 1$ also in this case.
\end{proof}

%
%

\begin{proof}[Proof of Theorem \ref{thm:str}]
Denote the Zhegalkin polynomial of $f$ by $\mathfrak{p}$. It is straightforward to verify that if $\mathfrak{p}$ has one of the special forms listed in the statement of the theorem, then $f$ does not have a variable identification minor of essential arity $\ess f - 1$ but it has one of essential arity $\ess f - 2$. For the converse implication, we will prove by induction on $\ess f$ that if $\mathfrak{p}$ is not of any of the special forms, then there is a variable identification minor $g$ of $f$ such that $\ess g = \ess f - 1$, i.e., $f$ has arity gap $1$.

If $\ess f = 2$ and $\mathfrak{p}$ is not of any of the special forms, then $\mathfrak{p} = x_i x_j + c$ or $\mathfrak{p} = x_i x_j + x_i + x_j + c$ where $c \in \{0,1\}$, and in both cases $\minor{\mathfrak{p}}{j}{i} = x_i + c$. In this case $\gap f = 1$.

If $\ess f = 3$, then $\mathfrak{p}$ has one of the following forms
\begin{align*}
& x_i x_j x_k + x_i x_j + x_i x_k + x_j x_k + a_i x_i + a_j x_j + a_k x_k + c, \\
& x_i x_j x_k + x_i x_k + x_j x_k + a_i x_i + a_j x_j + a_k x_k + c, \\
& x_i x_j x_k + x_i x_j + a_i x_i + a_j x_j + a_k x_k + c, \\
& x_i x_j + x_i x_k + x_j x_k + x_k + c, \\
& x_i x_j + x_i x_k + x_j x_k + x_i + x_j + x_k + c, \\
& x_i x_j + x_i x_k + a_i x_i + a_j x_j + a_k x_k + c, \\
& x_i x_k + a_i x_i + a_j x_j + a_k x_k + c,
\end{align*}
where $a_i, a_j, a_k, c \in \{0,1\}$. It is easy to verify that in each case $\minor{\mathfrak{p}}{j}{i}$ contains the term $x_i x_k$, and hence both $x_i$ and $x_k$ are essential in $\minor{f}{j}{i}$, and so $\gap f = 1$.

For the sake of induction, assume then that the claim holds for $2 \leq \ess f < n$, $n \geq 4$. Consider the case that $\ess f = n$. Since the case where $\deg \mathfrak{p} = 1$ is ruled out by the assumption that $\mathfrak{p}$ does not have any of the special forms and the case where $\deg \mathfrak{p} = 2$ is settled by Lemma \ref{lemma:deg2}, we can assume that $\deg \mathfrak{p} \geq 3$. Choose a variable $x_m$ from a term of the highest possible degree in $\mathfrak{p}$, and write
\[
\mathfrak{p} = x_m \mathfrak{q} + \mathfrak{r},
\]
where the polynomials $\mathfrak{q}$ and $\mathfrak{r}$ do not contain $x_k$. We clearly have that $\deg \mathfrak{q} = \deg \mathfrak{p} - 1$, and $\mathfrak{q}$ and $\mathfrak{r}$ represent functions with less than $n$ essential variables. Of course, every essential variable of $f$ except for $x_m$ occurs in $\mathfrak{q}$ or $\mathfrak{r}$. We have three different cases to consider, depending on the comparability under inclusion of the sets of variables occurring in $\mathfrak{q}$ and $\mathfrak{r}$.

\emph{Case 1.} Assume that there is a variable $x_i$ that occurs in $\mathfrak{q}$ but does not occur in $\mathfrak{r}$, and there is a variable $x_j$ that occurs in $\mathfrak{r}$ but does not occur in $\mathfrak{q}$. Write
\[
\mathfrak{q} = x_i \mathfrak{q}' + \mathfrak{q}'', \qquad \mathfrak{r} = x_j \mathfrak{r}' + \mathfrak{r}'',
\]
where $\mathfrak{q}'$, $\mathfrak{q}''$, $\mathfrak{r}'$, $\mathfrak{r}''$ do not contain $x_i$, $x_j$. Then
\[
\mathfrak{p} = x_m x_i \mathfrak{q}' + x_m \mathfrak{q}'' + x_j \mathfrak{r}' + \mathfrak{r}'',
\]
and we have that
\[
\minor{\mathfrak{p}}{j}{i} = x_m x_i \mathfrak{q}' + x_m \mathfrak{q}'' + x_i \mathfrak{r}' + \mathfrak{r}'',
\]
where no terms can cancel. Hence all essential variables of $f$ except for $x_j$ are essential in $\minor{f}{j}{i}$ and so $\gap f = 1$.

\emph{Case 2.} Assume that every variable occurring in $\mathfrak{r}$ occurs in $\mathfrak{q}$. In this case $\mathfrak{q}$ represents a function $q$ of essential arity $\ess f - 1$, containing all essential variables of $f$ except for $x_m$. We also have that $\deg \mathfrak{q} = \deg \mathfrak{p} - 1 \geq 2$.

\emph{Subcase 2.1.} If $\ess f \geq 5$, then $\ess q \geq 4$, and we can apply the inductive hypothesis, which tells us that there are variables $x_i$ and $x_j$ such that $\ess \minor{q}{i}{j} = \ess q - 1$. Hence $\minor{f}{i}{j}$ is represented by the polynomial $\minor{\mathfrak{p}}{i}{j} = x_m \minor{\mathfrak{q}}{i}{j} + \minor{\mathfrak{r}}{i}{j}$, and all essential variables of $f$ except for $x_i$ occur in $\minor{\mathfrak{p}}{i}{j}$, since no terms can cancel between $x_m \minor{\mathfrak{q}}{i}{j}$ and $\minor{\mathfrak{r}}{i}{j}$. Thus $\gap f = 1$.

\emph{Subcase 2.2.} If $\ess f = 4$, then $\ess q = 3$, and we can apply the inductive hypothesis as above unless $\mathfrak{q} = x_i x_j + x_i x_k + x_j x_k + c$ or $\mathfrak{q} = x_i x_j + x_i x_k + x_j x_k + x_i + x_j + c$. If this is the case, consider first the case where $\mathfrak{q}$ contains a variable $x_l \in \{x_i, x_j, x_k\}$ that does not occur in $\mathfrak{r}$. Consider then
\[
\minor{\mathfrak{p}}{m}{l} = x_l \mathfrak{q} + \mathfrak{r}.
\]
Then $x_l \mathfrak{q}$ contains the term $x_i x_j x_k$, which cannot be cancelled. Namely, all other terms of $x_l \mathfrak{q}$ have degree at most $2$, and since there are at most two variables occurring in $\mathfrak{r}$, the terms of $\mathfrak{r}$ also have degree at most $2$. Thus, all variables of $f$ except for $x_m$ occur in $\minor{\mathfrak{p}}{m}{l}$, and so the arity gap of $f$ is $1$.

Consider then the case that $\mathfrak{q}$ and $\mathfrak{r}$ contain the same variables, i.e., $x_i$, $x_j$, $x_k$. If $\deg \mathfrak{r} \leq 2$, then it is easily seen that $\minor{\mathfrak{p}}{m}{i}$ contains the term $x_i x_j x_k$, and all essential variables of $f$ except for $x_m$ are essential in $\minor{f}{m}{i}$. Otherwise, we can apply the inductive hypothesis on the function $r$ represented by $\mathfrak{r}$ and we obtain variables $x_\alpha$ and $x_\beta$ such that $\ess \minor{r}{\alpha}{\beta} = \ess r - 1$. It can be easily verified that no identification of variables brings $\mathfrak{q}$ into the zero polynomial, so $x_m$ and two other variables will occur in $\minor{\mathfrak{p}}{\alpha}{\beta} = x_m \minor{\mathfrak{q}}{\alpha}{\beta} + \minor{\mathfrak{r}}{\alpha}{\beta}$. We have that $\gap f = 1$ also in this case.

\emph{Case 3.} Assume that every variable occurring in $\mathfrak{q}$ occurs in $\mathfrak{r}$ but there is a variable $x_l$ that occurs in $\mathfrak{r}$ but does not occur in $\mathfrak{q}$. If $\deg \mathfrak{r} = 1$, then $\mathfrak{r} = x_l + \mathfrak{r}'$ where $\mathfrak{r}'$ does not contain $x_l$. Then $\minor{\mathfrak{p}}{m}{l} = x_l \mathfrak{q} + x_l + \mathfrak{r}'$, where the only term that may cancel out is $x_l$, and this happens if $\mathfrak{q}$ has a constant term $1$. Nevertheless, $x_l$ occurs in $\minor{\mathfrak{r}}{m}{l}$ because $\deg \mathfrak{q} \geq 2$. Of course, all other essential variables of $f$ except for $x_m$ also occur in $\minor{\mathfrak{p}}{m}{l}$, so $\gap f = 1$. We may thus assume that $\deg \mathfrak{r} \geq 2$.

\emph{Subcase 3.1.} Assume first that $\ess f = 4$ (in which case $\mathfrak{r}$ contains three variables and $\mathfrak{q}$ contains at most two variables) and $\mathfrak{r} = x_i x_j + x_i x_k + x_j x_k + c$ or $\mathfrak{r} = x_i x_j + x_i x_k + x_j x_k + x_i + x_j + c$. Since we assume that $\deg \mathfrak{p} \geq 3$, we have that $\deg \mathfrak{q} \geq 2$ and hence $\mathfrak{q}$ contains at least two variables. Thus exactly two variables occur in $\mathfrak{q}$ and so also $\deg \mathfrak{q} = 2$. Then $\mathfrak{q} = x_\alpha x_\beta + b_1 x_\alpha + b_2 x_\beta + d$ where $\alpha, \beta \in \{i,j,k\}$ and $b_1, b_2, d \in \{0,1\}$. Let $\gamma \in \{i, j, k\} \setminus \{\alpha, \beta\}$. Then $\minor{\mathfrak{p}}{m}{\gamma}$ contains the term $x_i x_j x_k$, and hence all essential variables of $f$ except for $x_m$ occur in $\minor{\mathfrak{p}}{m}{\gamma}$, and so $\gap f = 1$.

\emph{Subcase 3.2.} Assume then that $\ess f > 4$ or $\ess f = 4$ but $\mathfrak{r}$ does not have any of the special forms. In this case we can apply the inductive hypothesis on the function $r$ represented by $\mathfrak{r}$. Let $x_i$ and $x_j$ be such that $\ess \minor{r}{j}{i} = \ess r - 1$. If $\minor{\mathfrak{q}}{j}{i} \neq 0$, then $x_m$ and all other essential variables of $f$ except for $x_j$ occur in $\minor{\mathfrak{p}}{j}{i}$, and we are done---the arity gap of $f$ is $1$. We may thus assume that $\minor{\mathfrak{q}}{j}{i} = 0$. Write $\mathfrak{q}$ and $\mathfrak{r}$ in the form
\begin{align*}
\mathfrak{q} &= x_i x_j \mathfrak{a}_1 + x_i \mathfrak{a}_2 + x_j \mathfrak{a}_3 + \mathfrak{a}_4, \\
\mathfrak{r} &= x_i x_j \mathfrak{b}_1 + x_i \mathfrak{b}_2 + x_j \mathfrak{b}_3 + \mathfrak{b}_4,
\end{align*}
where the polynomials $\mathfrak{a}_1$, $\mathfrak{a}_2$, $\mathfrak{a}_3$, $\mathfrak{a}_4$, $\mathfrak{b}_1$, $\mathfrak{b}_2$, $\mathfrak{b}_3$, $\mathfrak{b}_4$ do not contain $x_i$, $x_j$. Define the polynomials $\mathfrak{q}_1$, \ldots, $\mathfrak{q}_7$ as follows (cf.\ the proof of Theorem 4 in Salomaa \cite{Salomaa}):

$\mathfrak{q}_1$ consists of the terms common to $\mathfrak{a}_1$, $\mathfrak{a}_2$, and $\mathfrak{a}_3$.

$\mathfrak{q}_i$, $i = 2, 3$, consists of those terms common to $\mathfrak{a}_1$ and $\mathfrak{a}_i$ which are not in $\mathfrak{q}_1$.

$\mathfrak{q}_4$ consists of those terms common to $\mathfrak{a}_2$ and $\mathfrak{a}_3$ which are not in $\mathfrak{q}_1$.

$\mathfrak{q}_{4+i}$, $i = 1, 2, 3$, consists of the remaining terms in $\mathfrak{a}_i$.

Define the polynomials and $\mathfrak{r}_1$, \ldots, $\mathfrak{r}_7$ similarly in terms of the $\mathfrak{b}_i$'s. Note that for any $i \neq j$, $\mathfrak{q}_i$ and $\mathfrak{q}_j$ do not have any terms in common, and similarly $\mathfrak{r}_i$ and $\mathfrak{r}_j$ do not have any terms in common. Hence,
\begin{align*}
\mathfrak{q} &= x_i x_j (\mathfrak{q}_1 + \mathfrak{q}_2 + \mathfrak{q}_3 + \mathfrak{q}_5) +
x_i (\mathfrak{q}_1 + \mathfrak{q}_2 + \mathfrak{q}_4 + \mathfrak{q}_6) +
x_j (\mathfrak{q}_1 + \mathfrak{q}_3 + \mathfrak{q}_4 + \mathfrak{q}_7) + \mathfrak{a}_4, \\
\mathfrak{r} &= x_i x_j (\mathfrak{r}_1 + \mathfrak{r}_2 + \mathfrak{r}_3 + \mathfrak{r}_5) +
x_i (\mathfrak{r}_1 + \mathfrak{r}_2 + \mathfrak{r}_4 + \mathfrak{r}_6) +
x_j (\mathfrak{r}_1 + \mathfrak{r}_3 + \mathfrak{r}_4 + \mathfrak{r}_7) + \mathfrak{b}_4.
\end{align*}
Identification of $x_i$ with $x_j$ yields
\begin{align*}
\minor{\mathfrak{q}}{j}{i} &= x_i (\mathfrak{q}_1 + \mathfrak{q}_5 + \mathfrak{q}_6 + \mathfrak{q}_7) + \mathfrak{a}_4, \\
\minor{\mathfrak{r}}{j}{i} &= x_i (\mathfrak{r}_1 + \mathfrak{r}_5 + \mathfrak{r}_6 + \mathfrak{r}_7) + \mathfrak{b}_4.
\end{align*}
Since we are assuming that $\minor{\mathfrak{q}}{j}{i} = 0$, we have that $\mathfrak{q}_1 = \mathfrak{q}_5 = \mathfrak{q}_6 = \mathfrak{q}_7 = \mathfrak{a}_4 = 0$. On the other hand, $\mathfrak{q} \neq 0$, so $\mathfrak{q}_2$, $\mathfrak{q}_3$, $\mathfrak{q}_4$ are not all zero. Thus
\[
\mathfrak{q} = x_i x_j (\mathfrak{q}_2 + \mathfrak{q}_3) + x_i (\mathfrak{q}_2 + \mathfrak{q}_4) + x_j (\mathfrak{q}_3 + \mathfrak{q}_4).
\]
All essential variables of $f$ except for $x_j$ are contained in $\minor{\mathfrak{r}}{j}{i}$.

\emph{Subcase 3.2.1.} Assume that there is a variable $x_t$ occurring in $\mathfrak{b}_4$ that does not occur in $\mathfrak{r}_1$, $\mathfrak{r}_5$, $\mathfrak{r}_6$, $\mathfrak{r}_7$. Consider
\[
\minor{\mathfrak{p}}{m}{t} = x_t \mathfrak{q} + \mathfrak{r} = x_l \mathfrak{q} + x_i x_j \mathfrak{b}_1 + x_i \mathfrak{b}_2 + x_j \mathfrak{b}_3 + \mathfrak{b}_4.
\]
Cancelling may only happen between a term of $x_t \mathfrak{q}$ and a term of $\mathfrak{r}$. No term of $\mathfrak{b}_4$ can be cancelled, because every term of $x_t \mathfrak{q}$ contains $x_i$ or $x_j$ but the terms of $\mathfrak{b}_4$ do not contain either. The variables that do not occur in $\mathfrak{b}_4$ occur in some terms of $\mathfrak{b}_1$, $\mathfrak{b}_2$, $\mathfrak{b}_3$ that do not contain $x_t$. Thus, all essential variables of $f$ except for $x_m$ occur in $\minor{\mathfrak{p}}{m}{t}$, and so in this case $f$ has arity gap $1$.

\emph{Subcase 3.2.2.} Assume that all variables of $\mathfrak{r}$ except for $x_i$, $x_j$ occur already in $\mathfrak{r}_1 + \mathfrak{r}_5 + \mathfrak{r}_6 + \mathfrak{r}_7$. Consider
\begin{equation}
\label{eq:case3a}
\begin{split}
\minor{\mathfrak{p}}{m}{i} = x_i x_j &(\mathfrak{q}_2 + \mathfrak{q}_4 + \mathfrak{r}_1 + \mathfrak{r}_2 + \mathfrak{r}_3 + \mathfrak{r}_5) + {} \\
x_i &(\mathfrak{q}_2 + \mathfrak{q}_4 + \mathfrak{r}_1 + \mathfrak{r}_2 + \mathfrak{r}_4 + \mathfrak{r}_6) + {} \\
x_j &(\mathfrak{r}_1 + \mathfrak{r}_3 + \mathfrak{r}_4 + \mathfrak{r}_7) + \mathfrak{b}_4.
\end{split}
\end{equation}

\emph{Subcase 3.2.2.1.} Assume first that $x_i$ does not occur in $\minor{\mathfrak{p}}{m}{i}$ in \eqref{eq:case3a}. Then
\begin{align*}
\mathfrak{q}_2 + \mathfrak{q}_4 + \mathfrak{r}_1 + \mathfrak{r}_2 + \mathfrak{r}_3 + \mathfrak{r}_5 = 0, \\
\mathfrak{q}_2 + \mathfrak{q}_4 + \mathfrak{r}_1 + \mathfrak{r}_2 + \mathfrak{r}_4 + \mathfrak{r}_6 = 0,
\end{align*}
and since the $\mathfrak{r}_i$'s do not have terms in common, we have that
\[
\mathfrak{r}_1 + \mathfrak{r}_2 = \mathfrak{q}_2 + \mathfrak{q}_4, \qquad \mathfrak{r}_3 = \mathfrak{r}_4 = \mathfrak{r}_5 = \mathfrak{r}_6 = 0.
\]
Then all variables of $\mathfrak{r}$ except for $x_i$, $x_j$ occur already in $\mathfrak{r}_1 + \mathfrak{r}_7$. Consider
\begin{equation}
\label{eq:case3b}
\begin{split}
\minor{\mathfrak{p}}{m}{j} = x_i x_j &(\mathfrak{q}_3 + \mathfrak{q}_4 + \mathfrak{r}_1 + \mathfrak{r}_2 + \mathfrak{r}_3 + \mathfrak{r}_5) + {} \\
x_i &(\mathfrak{r}_1 + \mathfrak{r}_2 + \mathfrak{r}_4 + \mathfrak{r}_6) + {} \\
x_j &(\mathfrak{q}_3 + \mathfrak{q}_4 + \mathfrak{r}_1 + \mathfrak{r}_3 + \mathfrak{r}_4 + \mathfrak{r}_7) + \mathfrak{b}_4 \\
{} = x_i x_j &(\mathfrak{q}_2 + \mathfrak{q}_3) + {} \\
x_i &(\mathfrak{r}_1 + \mathfrak{r}_2) + {} \\
x_j &(\mathfrak{q}_2 + \mathfrak{q}_3 + \mathfrak{r}_2 + \mathfrak{r}_7) + \mathfrak{b}_4.
\end{split}
\end{equation}
All variables of $\mathfrak{r}_1$ are there on the fifth line of \eqref{eq:case3b}. If a term of $\mathfrak{r}_7$ is cancelled by a term of $\mathfrak{q}_2 + \mathfrak{q}_3$ on the sixth line, it still remains on the fourth line, so all variables of $\mathfrak{r}_7$ are also there. We still need to verify that the variables $x_i$ and $x_j$ are not cancelled out from \eqref{eq:case3b}. If $\mathfrak{q}_2 + \mathfrak{q}_3 \neq 0$ then we are done. Assume then that $\mathfrak{q}_2 + \mathfrak{q}_3 = 0$, in which case $\mathfrak{q}_4 \neq 0$. Since
\[
\mathfrak{r}_1 + \mathfrak{r}_2 + \mathfrak{r}_4 + \mathfrak{r}_6 = \mathfrak{r}_1 + \mathfrak{r}_2 = \mathfrak{q}_2 + \mathfrak{q}_4 = \mathfrak{q}_4 \neq 0,
\]
we have $x_i$ in \eqref{eq:case3b}. Since
\[
\mathfrak{q}_3 + \mathfrak{q}_4 + \mathfrak{r}_1 + \mathfrak{r}_3 + \mathfrak{r}_4 + \mathfrak{r}_7 = \mathfrak{q}_4 + \mathfrak{r}_1 + \mathfrak{r}_7
\]
and $\mathfrak{r}_1 + \mathfrak{r}_7$ contain all variables of $\mathfrak{r}$ except for $x_i$, $x_j$, but $\mathfrak{q}_4$ does not, $ \mathfrak{q}_4 + \mathfrak{r}_1 + \mathfrak{r}_7 \neq 0$, so we also have $x_j$ in \eqref{eq:case3b}. Thus, the arity gap of $f$ equals $1$ in this case.

\emph{Subcase 3.2.2.2.} Assume then that $x_i$ occurs in $\minor{\mathfrak{p}}{m}{i}$ in \eqref{eq:case3a}. Nothing cancels out on the third line of \eqref{eq:case3a}, and therefore the variables of $\mathfrak{r}_1$ and $\mathfrak{r}_7$ occur in $\minor{\mathfrak{p}}{m}{i}$. Terms of $\mathfrak{r}_5$ may be cancelled out by terms of $\mathfrak{q}_2 + \mathfrak{q}_4$ on the first line of \eqref{eq:case3a} but such terms will remain on the second line. Thus the variables of $\mathfrak{r}_5$ occur in $\minor{\mathfrak{p}}{m}{i}$. A similar argument shows that the variables of $\mathfrak{r}_6$ also occur in $\minor{\mathfrak{p}}{m}{i}$. In order for $f$ to have arity gap $1$, we still need to verify that $x_j$ occurs in $\minor{\mathfrak{p}}{m}{i}$. If $\mathfrak{q}_2 + \mathfrak{q}_4 + \mathfrak{r}_1 + \mathfrak{r}_2 + \mathfrak{r}_3 + \mathfrak{r}_5 \neq 0$, then we are done. We may thus assume that
\begin{equation}
\label{eq:1stline}
\mathfrak{q}_2 + \mathfrak{q}_4 + \mathfrak{r}_1 + \mathfrak{r}_2 + \mathfrak{r}_3 + \mathfrak{r}_5 = 0.
\end{equation}
By the assumption that $x_i$ occurs in $\minor{\mathfrak{p}}{m}{i}$, the second line of \eqref{eq:case3a} does not vanish, i.e.,
\[
0 \neq \mathfrak{q}_2 + \mathfrak{q}_4 + \mathfrak{r}_1 + \mathfrak{r}_2 + \mathfrak{r}_4 + \mathfrak{r}_6 = \mathfrak{r}_3 + \mathfrak{r}_4 + \mathfrak{r}_5 + \mathfrak{r}_6.
\]
If the third line of \eqref{eq:case3a} does not vanish either, i.e., $\mathfrak{r}_1 + \mathfrak{r}_3 + \mathfrak{r}_4 + \mathfrak{r}_7 \neq 0$, then we have both $x_i$ and $x_j$ and we are done. We may thus assume that $\mathfrak{r}_1 + \mathfrak{r}_3 + \mathfrak{r}_4 + \mathfrak{r}_7 = 0$, i.e., $\mathfrak{r}_1 = \mathfrak{r}_3 = \mathfrak{r}_4 = \mathfrak{r}_7 = 0$. Then all variables of $\mathfrak{r}$ except for $x_i$, $x_j$ occur already in $\mathfrak{r}_5 + \mathfrak{r}_6$. Equation \eqref{eq:1stline} implies that $\mathfrak{r}_2 + \mathfrak{r}_5 = \mathfrak{q}_2 + \mathfrak{q}_4$. Consider
\begin{equation}
\label{eq:case3c}
\begin{split}
\minor{\mathfrak{p}}{m}{j} =
x_i x_j &(\mathfrak{q}_3 + \mathfrak{q}_4 + \mathfrak{r}_1 + \mathfrak{r}_2 + \mathfrak{r}_3 + \mathfrak{r}_5) + {} \\
x_i &(\mathfrak{r}_1 + \mathfrak{r}_2 + \mathfrak{r}_4 + \mathfrak{r}_6) + {} \\
x_j &(\mathfrak{q}_3 + \mathfrak{q}_4 + \mathfrak{r}_1 + \mathfrak{r}_3 + \mathfrak{r}_4 + \mathfrak{r}_7) + \mathfrak{b}_4 \\
= x_i x_j &(\mathfrak{q}_2 + \mathfrak{q}_3) + {} \\
x_i &(\mathfrak{q}_2 + \mathfrak{q}_4 + \mathfrak{r}_5 + \mathfrak{r}_6) + {} \\
x_j &(\mathfrak{q}_3 + \mathfrak{q}_4) + \mathfrak{b}_4.
\end{split}
\end{equation}
Assume first that $\mathfrak{q}_2 + \mathfrak{q}_3 = 0$, in which case $\mathfrak{q}_4 \neq 0$. If a term of $\mathfrak{r}_5 + \mathfrak{r}_6$ is cancelled by a term of $\mathfrak{q}_4$ on the second line of \eqref{eq:case3c}, it will still remain on the third line. Therefore we have in $\minor{\mathfrak{p}}{m}{j}$ all variables of $\mathfrak{r}$ except for $x_i$ and $x_j$. Since $\mathfrak{r}_5 + \mathfrak{r}_6$ contains all variables of $\mathfrak{r}$ except for $x_i$, $x_j$ but $\mathfrak{q}_2 + \mathfrak{q}_4 = \mathfrak{q}_4$ does not, the second line of \eqref{eq:case3c} does not vanish, and so we have $x_i$. We also have $x_j$ because $\mathfrak{q}_3 + \mathfrak{q}_4 = \mathfrak{q}_4 \neq 0$ on the third line. In this case $f$ has arity gap $1$.

Assume then that $\mathfrak{q}_2 + \mathfrak{q}_3 \neq 0$. Then the first line of \eqref{eq:case3c} does not vanish and both $x_i$ and $x_j$ occur in $\minor{\mathfrak{p}}{m}{j}$. If any term of $\mathfrak{r}_5 + \mathfrak{r}_6$ is cancelled by a term of $\mathfrak{q}_2$ on the second line of \eqref{eq:case3c}, it still remains on the first line, and if it is cancelled by a term of $\mathfrak{q}_4$, it remains on the third line. Thus all variables of $\mathfrak{r}$ occur in $\minor{\mathfrak{p}}{m}{j}$, and $f$ has arity gap $1$ again. This completes the proof of Theorem \ref{thm:str}.
\end{proof}

\section{Concluding remarks}

We do not know whether the upper bound on arity gap given by Theorem \ref{thm:gen} is sharp. For base sets $A$ with $k \geq 3$ elements, we do not know whether there exists an operation $f$ on $A$ with $\ess f \geq k + 1$ and $\gap f \geq 3$. We know that for all $k \geq 2$, there are operations on a $k$-element set $A$ with arity gap $2$. Consider for instance the quasi-linear functions of Burle \cite{Burle}. A function $f$ is \emph{quasi-linear} if it has the form
\[
f = g(h_1(x_1) \oplus h_2(x_2) \oplus \dots \oplus h_n(x_n)),
\]
where $h_1, \ldots, h_n : A \to \{0,1\}$, $g : \{0,1\} \to A$ are arbitrary mappings and $\oplus$ denotes addition modulo 2. It is easy to verify that if those $h_i$'s that are nonconstant coincide (and $g$ is not a constant map), then $f$ has arity gap $2$.

In general, if there is an operation $f$ on a $k$-element set $A$ with with $\gap f = m$, then there are operations of arity gap $m$ on all sets $B$ of at least $k$ elements. Namely, it is easy to see that any operation $g$ on $B$ of the form
\[
g = \phi(f(\gamma(x_1), \gamma(x_2), \ldots, \gamma(x_n))),
\]
where $\gamma : B \to A$ is surjective and $\phi : A \to B$ is injective, satisfies $\ess g = \ess f$ and $\gap g = \gap f$.

\end{document}